\documentclass[12pt]{article}

\usepackage{amsmath,amssymb,enumerate,bbm}
\newtheorem{theorem}{Theorem}[section]

\newtheorem{problem}[theorem]{Problem}

\numberwithin{equation}{section}
\newenvironment{proof}{\noindent\textbf{Proof\ }}{\hspace*{\fill}$\Box$\medskip}

 \makeatletter
 
 \newcommand{\Rmnum}[1]{\expandafter\@slowromancap\romannumeral #1@}
 \makeatother

\begin{document}
\title{The number of occurrences of a fixed spread among $n$ directions in
vector spaces over finite fields}\author{Le Anh Vinh\\
Mathematics Department\\
Harvard University\\
Cambridge, MA 02138, US\\
vinh@math.harvard.edu}\maketitle

\begin{abstract}
We study a finite analog of a problem of Erd\"os, Hickerson and Pach on the maximum number of occurrences of a fixed angle among $n$ directions in three-dimensional spaces. 
\end{abstract}

\section{Introduction}

Let $\mathbbm{F}_q$ denote the finite field with
$q$ elements where $q \gg 1$ is an odd prime power. For any $x, y \in
\mathbbm{F}_q^d$, the distance between $x, y$ is defined as $\|x - y\|= (x_1 -
y_1)^2 + \ldots + (x_d - y_d)^2$. Let $E \subset \mathbbm{F}_q^d$, $d
\geqslant 2$. Then the finite analog of the classical Erd\"os distance problem is
to determine the smallest possible cardinality of the set
\begin{equation}
  \Delta (E) =\{\|x - y\|: x, y \in E\},
\end{equation}
viewed as a subset of $\mathbbm{F}_q$. Bourgain, Katz and Tao (\cite{bourgain-katz-tao}), showed,
using intricate incidence geometry, that for every $\varepsilon > 0$, there
exists $\delta > 0$, such that if $E \in \mathbbm{F}_q^2$ and
$C_{\varepsilon}^1 q^{\varepsilon} \leqslant |E| \leqslant C_{\varepsilon}^2
q^{2 - \varepsilon}$, then $| \Delta (E) | \geqslant C_{\delta} |E|^{\frac{1}{2}
+ \delta}$ for some constants $C_{\varepsilon}^1, C_{\varepsilon}^2$ and
$C_{\delta}$. The relationship between $\varepsilon$ and $\delta$ in their
argument is difficult to determine. Going up to higher dimension using
arguments of Bourgain, Katz and Tao is quite subtle. Iosevich and Rudnev \cite{iosevich-rudnev}
establish the following result using Fourier analytic method.

\begin{theorem}(\cite{iosevich-rudnev})
  Let $E \subset \mathbbm{F}_q^d$ such that $|E| \gtrsim C q^{d / 2}$ for
  $C$ sufficiently large. Then
  \begin{equation}
    | \Delta (E) | \gtrsim \min \left\{ q, \frac{|E|}{q^{(d - 1) / 2}}
    \right\} .
  \end{equation}
\end{theorem}

Iosevich and his collaborators investigated several related results using this
method in a series of papers \cite{covert,hart-iosevich-solymosi,hart-iosevich,hart-iosevich-koh-rudnev,iosevich-koh,iosevich-rudnev,iosevich-senger}. Using graph theoretic method, the author reproved some of these resutls in \cite{vinh-ejc1,vinh-ejc2,vinh-ejc,vinh-fkw,vinh-dg}. The advantages of the graph theoretic method are twofold. First, we can reprove and sometimes improve several known results in vector spaces over finite fields. Second, our approach works transparently in the non-Euclidean setting.  In this note, we use the graph theoretic method to study a finite analog of a related problem of Erdos, Hickerson and Pach \cite{ehp}.

\begin{problem} (\cite{ehp})
  Give a good asymptotic bounds for the maximum number of occurrences of a
  fixed angle $\gamma$ among $n$ unit vectors in three-dimensional spaces.
\end{problem}

If $\gamma = \pi / 2$, the maximum number of orthogonal pairs is known to be
$\Theta (n^{4 / 3})$ as this problem is equivalent to bounding the number of
point-line incidences in the plane (see \cite{brass} for a detailed discussion). For any
other angle $\gamma \neq \pi / 2$, we are far from haveing good estimates for
the maximum number of occurrences of $\gamma$. The only known upper bound is
still $O (n^{4 / 3})$. the same as for orthogonal pairs. For the lower bound,
Swanepoel and Valtr \cite{sv} established the bound $\Omega (n \log n)$, improving
an earlier result of Erdos, Hickerson and Pach \cite{ehp}. It is, however, widely
believed that the $\Omega (n \log n)$ lower bound can be much improved. 

The purpose of this note is to study an analog of this problem in the
three-dimension space over finite fields. In vector spaces over finite fields,
however, the separation of lines is not measured by the transcendental notion
of angle. A remarkable approach of Wildberger \cite{norman1,norman2} by recasting metrical
geometry in a purely algebraic setting, eliminate the difficulties in defining
an angle by using instead the notion of spread - in Euclidean geometry the
square of the sine of the angle between two rays lying on those lines (the
notation of spread will be defined precisely in Section 2). Using this
notation, we now can state the main result of this note.

\begin{theorem}\label{mt1-man}
  Let $E$ be a set of unit vectors in $\mathbbm{F}_q^3$ with $q^{3/2} \ll  |E| \ll q^2$. For any $\gamma \in \mathbbm{F}_q$, let $f_{\gamma}(E)$ denote the number of occurrences of a fixed spread $\gamma$ among $E$. Then $f_{\gamma} (E) =
  \Theta (|E|^2 / q)$ if $1-\gamma$ is a square in $\mathbbm{F}_q$ and $f_{\gamma}(E)=0$ otherwise.
\end{theorem}

The rest of this note is organized as follows. In Section 2, we follow Wildberger's construction of affine and projective rational trigonometry to define the notions of quadrance and spread. We then define the main tool of our proof, the finite Poincar\'e graphs. Using these graphs, we give a proof of Theorem \ref{mt1-man} in Section 3.

\section{Quadrance, Spread and finite Poincar\'e graphs}

In this section, we follow Wildberger's construction of affine and projective rational trigonometry over finite fields. Interested readers can see \cite{norman1,norman2} for a detailed discussion. 

\subsection{Quadrance and Spread: affine rational geometry}
We work in a three-dimensional vector space over a field $F$, not of characteristic two. Elements of the vector space are called points or vectors (these two terms are equivalent and will be used interchangeably) and are denoted by $U, V, W$ and so on. The zero vector or point is denote $O$. The unique line $l$ through distinct points $U$ and $V$ is denoted $UV$. For a non-zero point $U$ the line $OU$ is denoted $[U]$. Fix a symmetric bilinear form and represent it by $U \cdot V$. In terms of this form, the line $UV$ is perpendicular to the line $WZ$ precisely when $(V-U) \cdot (Z-W) = 0$. A point $U$ is a null point or null vector when $U \cdot U = 0$. The origin $O$ is always a null point, and there are others as well. 

The distance (or so-called \textit{quadrance} in Wildberger's construction) between the points $U$ and $V$ is the number
\begin{equation}\label{quadrance}
Q(U,V) = (V-U) \cdot (V-U).
\end{equation}
The line $UV$ is a null line precisely when $Q(U,V) = 0$, or equivalently when it is perpendicular to itself. 

In Euclidean geometry, the separation of lines is traditionally measured by the transcendental notion of \textit{angle}. The difficulities in defining an angle precisely, and in extending the concept over an arbitrarily field, are eliminated in rational trigonometry by using instead the notion of \textit{spread} - in Euclidean geometry the square of the sine of the angle between two rays lying on those lines. Precisely, the \textit{spread} between the non-null lines $UW$ and $VZ$ is the number
\begin{equation}\label{spread}
s(UW,VZ) = 1 - \frac{((W-U)\cdot (Z-V))^2}{Q(U,W)Q(V,Z)}.
\end{equation}
This depends only on the two lines, not the choice of points lying on them. The spread between two non-null lines is $1$ precisely when they are perpendicular. Given a large set $E$ of unit vectors in $\mathbbm{F}_q^3$, our aim is to study the number of occurences of a fixed spread $\gamma \in \mathbbm{F}_q$ among $E$. 

\subsection{Finite Poincar\'e graphs: projective rational geometry}

Fix a three-dimensional vector space over a field with a symmetric bilinear form $U \cdot V$ as in the previous subsection. A line though the origion $O$ will now be called a projective point and denoted by a small letter such as $u$. The space of such projective points is called $n$ dimensional projective space. If $V$ is a non-zero vector in the vector space, then $v = [V]$ denote the projective point $OV$. A projective point is a null projective point when some non-zero null point lies on it. Two projective points $u = [U]$ and $v = [V]$ are perpendicular when they are perpendicular as lines. 

The \textit{projective quadrance} between the non-null projective points $u = [U]$ and $v = [V]$ is the number 
\begin{equation}\label{p-quadrance}
q(u,v) = 1 - \frac{(U\cdot V)^2}{(U \cdot U)(V \cdot V)}.
\end{equation}
This is the same as the spread $s(OU,OV)$, and has the value $1$ precisely when the projective points are perpendicular. 

The \textit{projective spread} between the intersecting projective lines $wu = [W,U]$ and $wv = [W,V]$ is defined to be the spread between these intersecting planes:

\begin{equation}\label{p-spread}
S(wu,wv) = 1 - \frac{\left(\left(U-\frac{U\cdot W}{W \cdot W}W\right)\cdot \left(V - \frac{V\cdot W}{W \cdot W}W\right)\right)^2}{\left(\left(U-\frac{U\cdot W}{W \cdot W}W\right)\cdot \left(U-\frac{U\cdot W}{W \cdot W}W\right)\right)\left(\left(V - \frac{V\cdot W}{W \cdot W}W\right) \cdot \left(V - \frac{V\cdot W}{W \cdot W}W\right)\right)}
\end{equation}

This approach is entirely algebraic and elementary which allows one to formulate two dimensional hyperbolic geometry as a projective theorey over a general field. Precisely, over the real numbers, the projective quadrance in the projective rational model is the negative of the square of the hyperbolic sine of the hyperbolic distance between the corresponding points in the Poincar\'e model, and the projective spread is the square of the sine of the angle between corresponding geodesics in the Poincar\'e model (see \cite{norman2}). 

Let $\Omega$ be the set of square-type non-isotropic $1$-dimensional subspaces of $\mathbbm{F}_q^3$ then $|\Omega| = q(q+1)/2$. For a fixed $\gamma \in \mathbbm{F}_q$, the \textit{finite Poincar\'e graph} $P_q(\gamma)$ has vertices as the points in $\Omega$ and edges between vertices $[Z],[W]$ if and only if $s(OZ,OW) =\gamma$. These graphs can be viewed as a companion of the well-known (and well studied) finite upper half plane graphs (see \cite{terras} for a survey on the finite upper half plane graphs). From the definition of the spread, the finite Poincar\'e graph $P_q(\gamma)$ is nonempty if and only if $1-\lambda$ is a square in $\mathbbm{F}_q$. 

We have the orthogonal group $O_3(\mathbbm{F}_q)$ acts transitively on $\Omega$, and yields a symmetric association scheme $\Psi(O_3(\mathbbm{F}_q),\Omega)$ of class $(q+1)/2$. The relations of $\Psi(O_3(\mathbbm{F}_q),\Omega)$ are given by
\begin{eqnarray*}
 R_1 & = & \{([U],[V]) \in \Omega \times \Omega \mid  (U+V) \cdot (U+V) = 0\},\\
 R_i & = & \{([U],[V]) \in \Omega \times \Omega \mid (U+V) \cdot (U+V) = 2 + 2 \nu^{- (i
  - 1)} \} \, (2 \leqslant i \leqslant (q - 1) / 2)\\
 R_{(q+1)/2} & = & \{([U], [V]) \in \Omega \times \Omega \cdot (U+V) \cdot (U+V) = 2\},
\end{eqnarray*}
where $\nu$ is a generator of the field $\mathbbm{F}_q$ and we assume $U\cdot U = 1$ for all $[U] \in \Omega$ (see \cite{bannai-hao-song}, Section 6). Note that $\Psi(O_3(\mathbbm{F}_q),\Omega)$ is isomorphic to the association scheme $PGL(2,q)/D_{2(q-1)}$ where $D_{2(q-1)}$ is a dihedral sugroup of order $2(q-1)$.  The graphs $(\Omega,R_i)$ are not Ramanujan in general, but fortunately, they are asymptotic Ramanujan for large $q$. The following theorem summaries the results from \cite{bannai-shimabukuro-tanaka}, Section 2 in a rough form.

\begin{theorem} (\cite{bannai-shimabukuro-tanaka})\label{bhs} The graphs $(\Omega,R_i)$ ($1\leq i \leq (q+1)/2$) are regular of valency $Cq(1+o(1))$. Let $\lambda$ be any eigenvalue of the graph $(\Omega,R_i)$ with $\lambda \neq$ valency of the graph then
\[|\lambda| \leq c(1+o(1))\sqrt{q},\]
for some $C, c > 0$ (In fact, we can show that $c = 1/2$). 
\end{theorem}

Theorem \ref{bhs} implies that the finite Poincar\'e graphs $P_q(\gamma)$ are asymptotic Ramanujan whenever $1 - \gamma$ is a square in $\mathbbm{F}_q$. Precisely, we have the following theorem.

\begin{theorem} \label{spectral-man}
a) If $1 - \gamma$ is not a square in $\mathbbm{F}_q$ then the finite Poincar\'e graph $P_q(\gamma)$ is empty. 

b) If $1 - \gamma$ is a square in $\mathbbm{F}_q$ then the finite Poincar\'e graph $P_q(\gamma)$ is regular of valency $Cq(1+o(1))$. Let $\lambda$ be any eigenvalue of the graph $P_q(\gamma)$ with $\lambda \neq$ valency of the graph then
\[|\lambda| \leq c(1+o(1))\sqrt{q},\]
for some $C, c > 0$.
\end{theorem}

\begin{proof}
a) Suppose that $[U], [V] \in \Omega$ and $s(OU,OV) = \gamma$ then 
\[1-\gamma = \frac{(U\cdot V)^2}{(U \cdot U)(V \cdot V)}.\] 
But $U, V$ are square-type so $1-\gamma$ is a square in $\mathbbm{F}_q$.

b) It is easy to see that the finite Poincar\'e graphs $P_q(1-\nu^{2-2i}) = (\Omega,R_i)$ for $1\leq i \leq (q-1)/2$ and $P_q(1) = (\Omega,R_{(q+1)/2})$. The theorem follows immediately from Theorem \ref{bhs}.
\end{proof}

\section{Proof of Theorem \ref{mt1-man}}

We call a graph $G = (V, E)$ $(n, d,
\lambda)$-regular if $G$ is a $d$-regular graph on $n$ vertices with the absolute value of each of its eigenvalues but the
largest one is at most $\lambda$. It is well-known that if $\lambda \ll d$ then a $(n,d,\lambda)$-regular graph behaves similarly as a random graph $G_{n,d/n}$. Presicely, we have the following result (see Corollary 9.2.5 and Corollary 9.2.6 in \cite{alon-spencer}).

\begin{theorem} \label{expander} (\cite{alon-spencer})
  Let $G$ be a $(n, d, \lambda)$-regular graph.  For every set of vertices $B$ of $G$, we have
  \begin{equation}\label{f2}
    |e (B) - \frac{d}{2 n} |B|^2 | \leqslant \frac{1}{2} \lambda |B|,
  \end{equation}
  where $e (B)$ is number of edges in the induced subgraph of $G$ on $B$.
\end{theorem}

Let $E$ be a set of $m$ unit vectors in $\mathbbm{F}_q^3$ then $E$ can be viewed as a subset of $\Omega$. The number of occurrences of a fixed spread $\gamma$ among $E$ can be realized as the number of edges in the induced subgraph of the finite Poincar\'e graph $P_q(\gamma)$ on the vertex set $E$. Thus, from Theorem \ref{spectral-man}, $f_{\gamma}(E) = 0$ if $1 - \gamma$ is not a square in $\mathbbm{F}_q$. 

Suppose that $1 - \gamma$ is a square in $\mathbbm{F}_q$.  From Theorem \ref{spectral-man} and Theorem \ref{expander}, we have

\begin{equation}\label{bound-man}
|f_{\gamma}(E) - \frac{Cq(1+o(1))}{q(q+1)/2}|E|^2| \leq \frac{1}{2} c(1+o(1))\sqrt{q}|E|.
\end{equation}

Since $|E| \gg q^{3/2}$, we have $\frac{1}{2} c(1+o(1))\sqrt{q}|E| \ll \frac{Cq(1+o(1))}{q(q+1)/2}|E|^2$ and the theorem follows.


\end{document}